\documentstyle[11pt]{article}

\title{
{\bf Simple Heuristics for Unit Disk Graphs}$\,$\footnote{
An extended abstract containing
some of the results in this paper appears in the Proceedings of the
4th Canadian Conference on Computational Geometry, St. Johns,
Newfoundland, Canada, August 1992, pp 244-249.}
}

\author{
{\em M.V. Marathe}$\,$\thanks{
Department of Computer Science, University at Albany -- SUNY,
Albany, NY 12222.
Email: {\tt madhav@cs.albany.edu}.
Research supported by NSF Grant CCR-89-03319.}
\hspace*{0.15in} 
{\em H. Breu}$\,$\thanks{
Department of Computer Science,
University of British Columbia,
Vancouver, British Columbia,
Canada V6T 1Z4.
Email: {\tt breu@cs.ubc.ca}.}
\hspace*{0.15in}
{\em H.B. Hunt III}$\,$\thanks{
Department of Computer Science, University at Albany -- SUNY,
Albany, NY 12222.
Email: {\tt hbh@cs.albany.edu}.
Research supported by NSF Grant CCR-89-03319.}
\hspace*{0.15in} 
{\em S.S. Ravi}$\,$\thanks{
Department of Computer Science, University at Albany -- SUNY,
Albany, NY 12222.
Email: {\tt ravi@cs.albany.edu}.
Research supported by NSF Grant CCR-89-05296.}
\hspace*{0.15in} 
{\em D.J. Rosenkrantz}$\,$\thanks{
Department of Computer Science, University at Albany -- SUNY,
Albany, NY 12222.
Email: {\tt djr@cs.albany.edu}.
Research supported by NSF Grant CCR-90-06396.}
}

\date{}

\begin{document}

\newcommand{\QED}{$\Box$}

\newcommand{\newthmwithin}[3]{\newtheorem{#1q}{#2}[#3]
                        \newenvironment{#1}{\begin{#1q}\sf}{\end{#1q}}}

\newcommand{\newthm}[3]{\newtheorem{#1q}[#2q]{#3}
                        \newenvironment{#1}{\begin{#1q}\sf}{\end{#1q}}}
\newcommand{\newthmm}[3]{\newtheorem{#1q}[#2q]{#3}
                        \newenvironment{#1}{\begin{#1q}\rm}}

\newtheorem{theorem}{Theorem}[section]
\newtheorem{lemma}{Lemma}[section]
\newtheorem{corollary}{Corollary}[section]
\newtheorem{fact}{Fact}[section]
\newcommand{\makenewheading}[1]{\begin{tabbing} {\bf #1:} \end{tabbing}}
\newcommand{\set}[1]{\{ #1 \}}

\newcommand{\eq}[1]{(\ref{#1})}
\newcommand{\cp}[1]{\mbox{\sc cap}(#1)}
\newcommand{\dem}[1]{\mbox{\sc dem}(#1)}

\newenvironment{proof}{{\bf Proof:}}{}
\newenvironment{sketchproof}{{\bf Sketch of Proof:}}{}

\newcommand{\order}[1]{$\mbox{O}\left(#1\right)$}                   

\newcommand{\peqnp}{${\bf P}={\bf NP}$}
\newcommand{\pneqnp}{${\bf P}\neq{\bf NP}$}
\newcommand{\pep}{\mbox{$P_{\epsilon}$}}
\newcommand{\vone}{$v_{1}$}
\newcommand{\vtwo}{$v_{2}$}
\newcommand{\vi}{$v_{i}$}
\newcommand{\vj}{$v_{j}$}
\newcommand{\vn}{$v_{n}$}
\newcommand{\xzero}{$x_{0}$}
\newcommand{\xone}{$x_{1}$}
\newcommand{\xtwo}{$x_{2}$}
\newcommand{\xrmo}{$x_{r-1}$}
\newcommand{\xsi}{$x_{i}$}
\newcommand{\pkb}{\mbox{{$P_{k}^{b}$}}}
\newcommand{\ppred}{\mbox{${\cal T}$}}

\newcommand{\xproj}{\mbox{${\ell_{\cal X}}$}}
\newcommand{\yproj}{\mbox{${\ell_{\cal Y}}$}}
\newcommand{\vproj}{\mbox{${\ell_{\cal V}}$}}
\newcommand{\wproj}{\mbox{${\ell_{\cal W}}$}}

\newcommand{\sxproj}{\mbox{${S_{\cal X}}$}}
\newcommand{\syproj}{\mbox{${S_{\cal Y}}$}}
\newcommand{\svproj}{\mbox{${S_{\cal V}}$}}
\newcommand{\swproj}{\mbox{${S_{\cal W}}$}}

\newcommand{\newspacing}{\baselineskip=1.4\normalbaselineskip}
\newcommand{\oldspacing}{\baselineskip=\normalbaselineskip}

\newspacing

\begin{titlepage}
\maketitle
\thispagestyle{empty}

\begin{center}
{\bf Abstract}
\end{center}

Unit disk graphs are intersection graphs of circles of unit
radius in the plane.
We present simple and provably good heuristics
for a number of classical NP-hard optimization problems on unit disk graphs.
The problems considered include 
maximum independent set, minimum vertex cover, minimum coloring and
minimum dominating set.
We also present an on-line coloring heuristic 
which achieves a competitive ratio of 6 for unit disk graphs.
Our heuristics do not need a geometric representation
of unit disk graphs.
Geometric representations are used only in establishing
the performance guarantees of the heuristics. 
Several of our approximation algorithms can be extended 
to intersection graphs of circles of arbitrary radii in the plane,
intersection graphs of regular polygons, and to intersection 
graphs of higher dimensional regular objects.

\end{titlepage}

\section{Introduction, motivation and summary of results}

Intersection graphs of geometric objects have been both widely studied 
and used to model many problems
in real life \cite{Ro78}.
In this paper we consider intersection graphs of regular 
polygons, emphasizing intersection graphs of unit disks. 
A graph is a 
{\bf unit disk graph} 
if and only if its vertices can be put in  one to
one correspondence with equisized circles in a plane in such a way that two
vertices are joined by an edge 
if and only if the corresponding circles intersect.
(It is assumed that tangent circles intersect.)
Without loss of generality, the radius of each circle (disk) is assumed
to be 1.
Unit disk graphs have been used to model broadcast networks
\cite{Ha80,Ka84,YWS84} and optimal facility location \cite{WK88}.
For example, 
the problem of placing $k$ facilities where proximity is undesirable
can be modeled as the problem of finding an independent set of size $k$
in a unit disk graph \cite{WK88}.
The problem of assigning distinct
frequencies for transmitters with intersecting ranges
corresponds to the minimum coloring problem for unit disk graphs.
The minimum dominating set problem corresponds to selecting
a minimum number of transmitters so that all the other stations 
are within the range of at least one of the chosen transmitters \cite{Ha80}.
As observed in \cite{CCJ90}, unit disk graphs are not perfect.
(An odd cycle of length five or more is a unit disk graph but not
perfect.)
Also, unit disk graphs are not planar.
(A clique of size five or more is a unit disk graph but not
planar.) 
Thus, in general, many of the known efficient algorithms 
for planar graphs and
perfect graphs may not be applicable to unit disk graphs.

All of the problems mentioned above and many others 
remain NP-hard for unit disk graphs \cite{CCJ90}.
Motivated by the practical importance of these problems, 
we present simple heuristics 
with provably good performance guarantees for a number
of such problems including the three mentioned above.
Our heuristics do not need a geometric
representation of a unit disk graph as part of the input.
Geometric representations are used only in establishing
some graph theoretic properties of unit disk graphs.
These properties, in turn, are used in deriving the
performance guarantees provided by the heuristics.
Similar properties hold for 
the intersection graphs
of other regular polygons and geometric
objects in higher dimensions. 
Consequently,  our heuristics
can be extended to those intersection graphs as well.

Our results are summarized below.
To provide a proper perspective, 
known results regarding performance guarantees for 
general graphs are also included.
\begin{enumerate}
\item
We present a heuristic with a performance guarantee of 3/2
for the minimum vertex
cover problem for unit disk graphs.
(The best known heuristic for the vertex cover problem for general graphs
provides a performance guarantee of 2 \cite{GJ79,BE85}.
It is also known that there is no 
polynomial time approximation scheme (PTAS) for the vertex cover problem
for general graphs,
unless \peqnp \cite{ALM+92}.)

\item
We present a simple heuristic with a performance
guarantee of 3 for
the maximum independent set problem for unit disk graphs.
(The maximum independent set for general graphs is notoriously
hard to approximate; unless \peqnp, 
there is an $\epsilon >0$ such that no polynomial time algorithm 
for the problem can provide a 
performance guarantee of $O(n^{\epsilon})$ \cite{ALM+92}.)

\item
We show that the off-line minimum vertex coloring problem can be 
approximated to within a factor of 3 
of the optimal value for unit disk graphs.
We note that unless \pneqnp, 
there is no PTAS for the minimum coloring problem for unit disk graphs
since the 3-coloring problem for unit disk graphs is NP-hard \cite{CCJ90}. 
We also give an {\em on-line} \cite{Ira90} coloring heuristic 
with a competitive ratio of 6 for unit disk graphs.
(For general graphs, unless \peqnp,
there is an $\epsilon > 0$ such that 
no polynomial time algorithm for the off-line 
minimum coloring problem can provide a 
performance guarantee of $O(n^{\epsilon})$ \cite{LY93}.
On-line coloring algorithms with constant
competitive ratios have been obtained previously for several other
classes of graphs \cite{HMR+94,Sl94,MHR93,Ki91,GL88,Sl87,KT81}.
It is also known that no on-line algorithm can provide a
competitive ratio of $o(\log{n})$ even for trees \cite{Ira90}.)

\item
We observe that a very simple heuristic 
provides a performance
guarantee of 5 for the minimum dominating set problem 
for unit disk graphs.
This heuristic 
also provides a performance guarantee of 5 
for the minimum independent domination problem.
The heuristic can be modified to obtain another
heuristic which provides a 
performance guarantee of 10 for both the minimum total domination
and the minimum connected domination problems.
(For general graphs, the minimum dominating set problem can be
approximated to within a factor $O(\log{n})$ of the optimal value
\cite{Jo74};
it is also known that no polynomial time algorithm
can provide a performance guarantee of $o(\log{n})$,
unless every problem in {\bf NP} can be solved
in deterministic time $O(n^{{\rm poly}\:\log{n}})$ \cite{LY93}.
For general graphs, it is also known that there is an $\epsilon > 0$
such that the minimum 
independent domination problem cannot
be approximated to within a factor of $O(n^{\epsilon})$
unless \peqnp \cite{Kan93,Ir91}.)

\item
For intersection graphs of circles of arbitrary radii in the plane,
we present heuristics with performance guarantees of 5/3, 6 and 5
respectively for vertex cover, off-line coloring and independent
problems respectively.
These heuristics are obtained as extensions of the corresponding
heuristics for unit disk graphs.
We also show how the heuristics can be extended to intersection
graphs of regular polygons.
\end{enumerate}

The remainder of the paper is organized as follows.
Section 2 contains the necessary definitions.
In Section 3, we establish several properties of unit disk graphs.
The approximation algorithms of Section 4 utilize these properties.
Section 5 contains extensions of our results to other
intersection graphs and Section 6 contains some concluding remarks. 

\section{Definitions }

We defined unit disk graphs as the intersection graphs of
sets of unit disks in the plane.
Throughout this paper, it is assumed that tangent circles intersect.
The model described above will be referred to as the
{\bf intersection model} \cite{CCJ90}.
Unit disk graphs can also be defined using
the {\bf proximity model} \cite{CCJ90} in which
the nodes of the graph are in one-to-one correspondence with
a set of points in the plane, and
two vertices are joined by an edge if and only if the distance between the
corresponding points is at most some specified bound. 
When considering a geometric representation, we often do not distinguish
between a vertex and its corresponding circle or point.

In the intersection model, we assume that the radius of 
each disk is 1.
It is easily seen that two such disks in the plane 
intersect if and only if the distance between their centers is at most 2.
Thus it is easy to translate a description 
of a given unit disk graph in the intersection
model into a description in the proximity model and vice versa
in linear time.
The recognition problem for unit disk graphs was
shown to be NP-hard in \cite{BK93}.
As already mentioned, none of our heuristics requires a geometric 
representation of a unit disk graph as part of the input.

An approximation algorithm for an optimization problem $\Pi$
provides a {\bf performance guarantee} of $\rho$ if for every
instance $I$ of $\Pi$, the solution value returned by the approximation
algorithm is within a factor $\rho$ of the optimal value for $I$.
A {\bf polynomial time approximation scheme} (PTAS) for problem $\Pi$
is a polynomial time approximation algorithm which given an instance
$I$ of $\Pi$ and an $\epsilon > 0$, returns a solution which is 
within a factor $(1+\epsilon)$ of the optimal value for $I$.

For the sake of completeness, we now define several graph 
theoretic parameters.
Given an undirected graph $G(V,E)$,
a {\bf minimum vertex cover} for $G$ is a smallest
cardinality subset $V'$ of $V$ 
such that for each edge $(x,y)$ in $E$,
at least one of $x$ and $y$ is in $V'$.
A {\bf maximum independent set} is a maximum cardinality subset $V'$ of $V$
such that there is no edge between any two vertices in $V'$.
A {\bf dominating set } is a subset $V'$ of $V$ such that each vertex
in $(V-V')$ has at least one neighbor in $V'$.
An {\bf independent dominating set} is a set of nodes that is
both an independent set and a dominating set.
A subset $V'$ of $V$ is a {\bf total dominating set}
if every node in $V$ has a neighbor in $V'$
(i.e., in addition to the nodes in $V-V'$, each node in $V'$ must
also be dominated by another node in $V'$).
A {\bf connected dominating set} is a dominating set $V'$ in which the vertex
induced subgraph on $V'$ is connected.

We end this section with some graph theoretic notation.
We use $K_{p,q}$ to denote the complete bipartite graph
with $p$ and $q$ nodes in the two sets of the bipartition.
We also use $K_r$ to denote the clique with $r$ nodes.
Given a graph $G(V,E)$ and a node $v$, we use $N(v)$
to denote the set of nodes adjacent to $v$; 
we refer to $N(v)$ as the {\bf neighborhood} of $v$.
For a graph $G(V,E)$ and a subset
$V' \subseteq V$, we use $G(V')$ to denote 
the subgraph of $G$ induced on $V'$.
For the remainder of this paper, we assume that $|V| = n$
and $|E| = m$.

\section{Some properties of unit disk graphs}

Many of our heuristics are based on a 
forbidden subgraph property of unit disk
graphs.
The proof of this property relies on a geometric observation
concerning packing of unit disks in the plane.  
(The problems of Packing and Covering have been of interest
to researchers for quite some time. See \cite{CS88} for more on this 
subject.)

\begin{lemma}\label{lem:pack}
Let $C$ be a circle of radius $r$ and 
and let $S$ be a set circles of radius $r$ such that 
every circle in $S$ intersects $C$ and no two circles in $S$
intersect each other.
Then, $|S| \leq 5$.
\end{lemma}

\noindent
{\bf Proof:} 
Suppose $|S| \geq 6$. 
Let $s_i$, $1 \leq i \leq 6$, denote the centers of
any six circles in $S$.
Let $c$ denote the center of $C$.
Denote the ray $\overrightarrow{c\,s_i}$ by $r_i$ ($1 \leq i \leq 6$).
Since there are six rays emanating from $c$,
there must at least one pair of rays $r_j$ and $r_k$
such that the angle between them is at most $60^{\circ}$.
Now, it can be verified that the distance between
$s_j$ and $s_k$ is at most $2r$, which implies that
circles centered at $s_j$ and $s_k$ intersect,
contradicting our assumption.
Thus $|S| \leq 5$.  \QED

An immediate consequence of Lemma \ref{lem:pack} is the following:

\begin{lemma}\label{lem:k16}
Let $G(V,E)$ be a unit disk graph.
Then $G$ cannot contain an induced subgraph isomorphic to $K_{1,6}$. \QED
\end{lemma}

Lemma \ref{lem:k16} indicates that in any unit disk graph,
the size of a maximum independent set in the 
subgraph of $G$ induced on the neighborhood
of any vertex is at most 5.
The following lemma, which can be proven
in a manner similar to that of Lemma \ref{lem:pack}, 
points out that the neighborhoods of certain nodes 
in a unit disk graph have even smaller independent sets.

\begin{lemma}\label{lem:semi}

Let $G$ be a unit disk graph, and let $v$ be a vertex
such that the unit disk corresponding to $v$
(in some model for $G$) has the smallest $X$-coordinate. 
The size of a maximum independent set in 
$G(N(v))$ is at most 3. \QED

\end{lemma}

Finally, the following property of unit disk graphs is useful in
establishing the performance of our heuristic for on-line vertex coloring.

\begin{lemma}\label{lem:delta_by_6}
A unit disk graph $G$ with maximum node degree $\Delta$ contains
a clique of size at least $\lceil \Delta/6 \rceil + 1$.
\end{lemma}

\noindent
{\bf Proof:} Consider a geometric representation of $G$ in which
each node of $G$ corresponds to a circle of radius 1.
Let $v$ be a node of maximum degree $\Delta$.
The centers of the circles corresponding to the $\Delta$
neighbors of $v$ are all within the circle $C$ of radius 2
with center at the center of the circle corresponding to $v$.
Thus at least $\lceil \Delta/6 \rceil$ of the neighbors of $v$
must lie within some $60^\circ$ sector of the circle $C$.
Since the maximum distance between a pair of points in that
segment is 2, it follows that the nodes corresponding to 
the points in that sector must be pairwise adjacent.
Thus, these $\lceil \Delta/6 \rceil$ nodes, together with
the node $v$, must form a clique of size  
$\lceil \Delta/6 \rceil + 1$.  \QED

\section{Approximation Algorithms}
We now proceed to give approximations to the various graph theoretic problems.
We start with a heuristic for the minimum vertex cover problem. 

\subsection{Minimum Vertex Cover}

The heuristic given here is 
essentially the same as the one given in
in Bar-Yehuda and Even \cite{BE82} for planar graphs but the analysis 
which leads to a 
performance guarantee of 1.5 is  different. 
It is interesting that an 
algorithm which provides a performance guarantee of 1.5 for planar graphs
provides the same performance guarantee 
for unit disk graphs as well.

Before presenting the details of the heuristic, we
mention some results which are used in the design and
analysis of the heuristic.
The first result deals with the {\bf Nemhauser-Trotter
decomposition} (NT decomposition) of a graph. 
The properties of this decomposition are stated
in the following result from \cite{NT75,BE85}.

\begin{theorem}\label{th:nt} \cite{NT75,BE85}
Given an undirected graph $G(V,E)$, there are two disjoint subsets
$P$ and $Q$ of $V$ 
such that the following properties hold.
\begin{enumerate}
\item
There exists an optimum cover $VC^*(G)$ for $G$
such that  $P \subseteq VC^*(G)$. 
\item
If $D$ is a vertex cover for $G(Q)$, 
then $D \cup P$ is a vertex cover for $G$.
\item
For any optimum vertex cover $VC^*(G)$ of $G$,
$|VC^*(G)| \geq  |P| + |Q|/2$.
\end{enumerate} 
Moreover, the sets $P$ and $Q$ can be found 
in polynomial time. \QED
\end{theorem}

The second result that is used in the heuristic is a theorem
due to Hochbaum \cite{Ho83}.
This theorem points out how a near-optimal
vertex cover can be obtained in some cases starting from an NT decomposition.
We have included the proof from \cite{Ho83} 
since our heuristic uses the method 
given in the proof.

\begin{theorem}\cite{Ho83} \label{th:hoch}
Let $G(V,E)$ be a node weighted graph and let $P$ and $Q$
denote the subsets obtained using an NT decomposition of $G$.
If $G(Q)$ can be colored using $k$ colors, then we can find
a vertex cover whose cardinality is at most $2(1- 1/k)$ times
that of an optimal vertex cover. 
\end{theorem}

\noindent
{\bf Proof:} Color the nodes of $G(Q)$ using $k$ colors.
This procedure partitions $Q$ into at most $k$ color classes.
Let $S$ be a color class of largest cardinality.
Clearly, $|S| \geq |Q|/k$.
It is easy to see from Theorem~\ref{th:nt} that 
the set $C = P \cup (Q - S)$ 
is a vertex cover for $G$.
Further, $|C| \leq |P| + (k-1)|Q|/k$
$= \displaystyle{\frac{2(k-1)}{k} 
       \left[\frac{k|P|}{2(k-1)} + \frac{|Q|}{2} \right]}$. 
By Theorem~\ref{th:nt}, any optimal vertex cover for $G$ has cardinality at 
least $|P| + |Q|/2$.
Thus, $C$ is a vertex cover for $G$
whose cardinality is no more than
$2(1-1/k)$ times that of an optimal vertex cover. \QED

Finally, the heuristic uses the following property of
triangle-free unit disk graphs.
(A graph is {\bf triangle-free} if it does not contain
a subgraph isomorphic to $K_3$.)  

\begin{lemma}\label{lem:triangle_free}
Any triangle-free unit disk graph can be colored using
4 colors.
\end{lemma}

\noindent
{\bf Proof:} 
We first observe that every triangle-free unit disk graph
has a node with degree at most 3.
To see this, let $G$ be a triangle-free unit disk graph.
Since $G$ is triangle-free, the neighborhood of each node
is an independent set.
Thus, by Lemma~\ref{lem:semi}, $G$ must have a node of degree
at most 3.

From the above observation and the fact that 
any induced subgraph
of a triangle-free
unit disk graph is also 
a triangle-free unit disk graph, it is easy to
see that the following is an efficient algorithm for 4-coloring such graphs: 
Given a triangle-free unit disk graph with $n$ nodes,
we find and remove a node $v$ of degree at most 3, recursively 4-color the
resulting graph containing $n-1$ nodes and then
assign a color to $v$ that has not been used for any of the
at most three neighbors of $v$. \QED

We note that the algorithm presented in the proof 
of Lemma~\ref{lem:triangle_free} does not
require a geometric representation of the unit disk graph.

\begin{figure}
\oldspacing
\rule{6.5in}{0.01in}

  \begin{enumerate}
      \item $V_1 = \emptyset ; V'= V$. 
      \item {\bf while} $G(V')$ contains triangles {\bf do} 
      \item {\bf begin}
           \begin{enumerate}
                \item Let $ X \subseteq V'$ be such that $G(X)$ is a triangle.
                \item $ V_1 = V_1 \cup X$.
                \item $ V' = V'-X$.
           \end{enumerate}
           \item {\bf end} 
       \item
            Obtain the sets $P$ and $Q$ by applying 
            Nemhauser-Trotter's algorithm on $G(V')$. 
           
       \item 
           Color the vertices of $G(Q)$ using 4 colors.
           Let $S$ denote a color class containing the largest number
           of nodes. 
       \item
           {\bf output} $VC(G) = V_1 \cup P \cup (Q-S)$ as the approximate 
           vertex cover for $G(V,E)$.
   \end{enumerate}
\caption{Details of Heuristic VCover}
\label{fig:heur_vc}
\rule{6.5in}{0.01in}
\newspacing
\end{figure}

We are now ready to discuss the heuristic for vertex cover.
The details of this heuristic, which we call VCover,
are given in Figure \ref{fig:heur_vc}.
The heuristic begins by removing all
triangles from the graph $G$. 
Every time a triangle is
found, we add all the three vertices of the triangle 
to the vertex cover ($V_1$) obtained so far. 
To the resulting triangle-free graph,
we apply NT decomposition
and obtain the sets $P$ and $Q$. 
In the next step we 4-color the graph $G(Q)$.
(This is possible because of Lemma \ref{lem:triangle_free}.)
We then remove the vertices in a color class $S$ of maximum cardinality
and return the set $V_1 \cup P \cup (Q-S)$  
as the approximate vertex cover.
To determine the performance guarantee of this heuristic, we
need the following result from \cite{BE85}.

\begin{lemma} \label{lem:local} \cite{BE85}
Let $r_1$ and $r_2$ denote the local ratios of two heur¸DÿÜïeÿª¸lcover problem. 
If $G_1$ denotes the graph obtained after applying the first
heuristic to a graph $G$,
then the performance of the algorithm which applies the two 
heuristics in succession is at most $\,\max \{ r_1(G), r_2(G_1) \}$. \QED
\end{lemma}

\begin{theorem}\label{th:perf_vc}
Let $G$ be a unit disk graph.
Let $VC^*(G)$ denote an optimal vertex cover for $G$
and let $VC(G)$ denote the vertex cover produced by
Heuristic VCover.
Then $|VC(G)| \leq 1.5 |VC^*(G)|$.
\end{theorem}

\noindent
{\bf Proof:}
Let us think of Heuristic VCover (Figure \ref{fig:heur_vc})
as operating in two phases,
where the first phase consists of Steps 1 through 4 and
the second phase consists of the remaining steps.
The first phase of the algorithm guarantees a local ratio of 1.5, because
any optimal algorithm needs to pick at least two vertices from each triangle
and the heuristic picked all three. 
From Theorem~\ref{th:nt}, 
Lemma~\ref{lem:triangle_free}, and Theorem~\ref{th:hoch}, it follows that
the second phase of the algorithm also guarantees a local ratio 
of 1.5.
The result now follows from Lemma \ref{lem:local}. \QED

The time complexity of Heuristic VCover is dominated by the time it takes to 
execute the Nemhauser-Trotter algorithm.

\subsection {Minimum Vertex Coloring}

\subsubsection{Off-line Coloring}
We present below a heuristic for coloring the vertices of a unit
disk graph such that the number of colors used is no more than
3 times the number of colors used by an optimal algorithm.
We use the following result due to Hochbaum \cite{Ho83}.

\begin{theorem}\cite{Ho83} \label{th:hoch_coloring}
Let $G =(V,E)$ be a graph and let $\delta(G)$ denote the largest
$\delta$ such that $G$ contains a subgraph in which every vertex has a 
degree at least $\delta$. Then $G$ can be colored using $\delta(G) +1$
colors. 
Moreover, such a coloring can be obtained in $O(|V|+|E|)$ time. \QED
\end{theorem}

Our heuristic for off-line coloring simply colors the given unit disk graph
$G$ with $\delta(G)+1$ colors using the algorithm mentioned in the above
theorem.
We now prove the performance guarantee of this heuristic.

\newcommand{\offcol}{\mbox{\em OFFCOL}}
\newcommand{\opt}{\mbox{\em OPT}}

\begin{theorem}\label{th:perf_off_color}
Let $G$ be a unit disk graph. 
Let $\offcol(G)$ and
$\opt(G)$ denote  the number of colors used by the above heuristic
and that used by an optimal algorithm respectively. Then
$\offcol(G) \leq 3 \; \opt(G). $
\end{theorem}

\noindent
{\bf Proof:}
Let $\delta(G)$ be as defined in Theorem~\ref{th:hoch_coloring}. 
Thus,
\begin{equation}
\offcol(G)  = \delta(G) +1. \label{eq:off_col1}
\end{equation}
We now prove the following claim which relates $\opt(G)$ and $\delta(G)$.

\noindent
{\bf Claim 1:} $\opt(G)  \geq \delta(G)/3 + 1$. 

\noindent
{\bf Proof:} Let $H$ be a subgraph of $G$ such that every vertex in
$H$ has a degree at least $\delta(G)$. 
Consider a unit disk representation of $H$, and let
$v^* \in H$ be a vertex such that the center of the unit disk 
corresponding to $v^*$ has the leftmost $X$ coordinate. 
Let $N_H(v^*)$ denote the neighborhood of $v^*$ in $H$. 
By our choice of $H$,
\begin{equation}
|N_H(v^*)| \geq \delta(G). \label{eq:off_col2}
\end{equation}
Recall from Lemma~\ref{lem:semi} that
the subgraph induced on the nodes in $N_H(v^*)$ has an independent
set of size at most 3. 
Consider any valid coloring of the graph $G$. 
Since all the vertices colored with the same color form an
independent set, in any valid coloring of the
subgraph induced on the set of nodes in $N_H(v^*)$, no more than
3 vertices can belong to the same color class.  
Therefore, any valid coloring of the subgraph induced by the set 
$N_H(v^*) \cup \{v^* \}$ must use at least $|N_H(v^*)|/3 +1$ colors.
In other words,
$\opt(G)$  $\geq$  $|N_H(v^*)|/3 +1$.
Now, from Equation~(\ref{eq:off_col2}), it follows that
\begin{equation}
\opt(G) \geq \delta(G)/3 + 1. \label{eq:off_col3}
\end{equation}
The theorem now follows from Equations~(\ref{eq:off_col1}) 
and (\ref{eq:off_col3}). \QED

\subsubsection{On-Line Coloring}
An instance of an off-line coloring problem consists of 
all the nodes and edges of a graph.
In the {\em on-line} version of the coloring problem,
vertices of a graph are presented one at a time.
When a vertex $v$ is presented, all the edges
from $v$ to the vertices presented earlier
are given.
An on-line coloring algorithm assigns a color to the current
vertex before the next vertex is presented.
The color assigned to the vertex must be different from
the colors assigned to its neighbors presented earlier.
Once a color is assigned to a vertex, the algorithm is not
allowed to change the color at a future time.
The performance guarantee provided by a heuristic is
defined as the worst-case ratio of the number of colors
used by the heuristic to the number of colors needed
in an optimal off-line coloring.
In the context of on-line problems,
the performance guarantee provided by a heuristic is
referred to as its {\bf competitive ratio}.

In this section, we observe that given a unit disk graph $G(V,E)$, 
there is an on-line coloring heuristic which achieves a 
competitive ratio of 6. 
The heuristic simply colors the vertices using the following
greedy strategy:

\begin{description}
\item[Heuristic Greedy:] Color each vertex with the first available color.
\end{description}

The following lemma 
provides an upper bound on the
number of colors used by the greedy strategy for an arbitrary graph.
The lemma can be proved easily by induction on the number of vertices.

\begin{lemma} \label{lem:oncol_greedy}
Let $G$ be an arbitrary graph with maximum node degree $\Delta(G)$.
Let $Greedy(G)$ denote the number of colors used by Greedy
to color $G$.
Then $Greedy(G) \leq \Delta(G) + 1$. \QED
\end{lemma}

The following theorem is an easy consequence of Lemmas 
\ref{lem:oncol_greedy} and \ref{lem:delta_by_6}.
 
\begin{theorem}\label{th:oncol}
Heuristic Greedy achieves a competitive ratio of at most 6. \QED
\end{theorem}

\subsection{Maximum Independent Set}

Since unit disk graphs are $K_{1,6}$ free, 
one can use the algorithm in \cite{Ho83}
to obtain a performance guarantee of 5. 
However, we can do better for unit disk graphs 
because of the additional 
geometric structure they possess.

Our heuristic, which we call IS,
provides a performance guarantee of 3.
This heuristic is based on Lemma~\ref{lem:semi}
which points out that every unit disk graph $G$ has
a node $v$ such that size of any independent set in
the subgraph induced on the
neighborhood $N(v)$ of $v$ is at most 3.
Given the nodes and edges of a unit disk graph,
such a node can be found in polynomial time.
(Obviously, $O(n^4)$ time suffices.)
This observation in conjunction with the fact that
any induced subgraph of a unit disk graph is also
a unit disk graph leads to the following simple
heuristic for finding a near-optimal independent set
in a unit disk graph.
The heuristic repeatedly finds a node $v$ such that the
maximum independent set in the subgraph induced on $N(v)$ is
at most 3, adds $v$ to the
independent set, and deletes $v$ and $N(v)$ from the graph.
These steps are repeated until the graph becomes empty.
We now show that
heuristic IS provides a performance guarantee of 3.

\newcommand{\IS}{\mbox{\em IS}}
\newcommand{\OPT}{\mbox{\em OPT}}

\begin{theorem}
Let $G$ be a unit disk graph, and let
$\IS(G)$ and $\OPT(G)$ denote respectively the approximate independent set
produced by {\rm IS}, and an optimal independent set for $G$. 
Then $|\IS(G)| \geq  |\OPT(G)|/3$. \QED
\end{theorem}

\noindent
{\bf Proof:} 
Define the (closed) {\em r-neighborhood} of a vertex $v$ in $\IS(G)$
to be the set $N(v) \cup \{ v \}$, where $N(v)$ is the neighborhood of
$v$ when the algorithm adds $v$ to $\IS(G)$.
By construction, every vertex in the graph, and therefore in $\OPT(G)$,
is in the r-neighborhood of at least one vertex in $\IS(G)$.
Also by construction, the size of a maximum independent set in every
r-neighborhood is 3.
Therefore, since $\OPT(G)$ is independent, every
r-neighborhood contains at most 3 vertices from $\OPT(G)$.
It follows that $|\OPT(G)| \leq 3 |\IS(G)|$. \QED

The running time of the heuristic can be reduced significantly
when a geometric representation of the unit disk graph is
available.
Recall from the proof of Lemma~\ref{lem:semi}
that a node $v$ corresponding to a unit disk whose center
has the least $X$-coordinate satisfies the required property.
Therefore, 
when a geometric representation of the unit disk graph is given,
we can implement the above heuristic in the following manner.
Form a list
$L$ of the centers of the unit disks in increasing order of their 
$X$-coordinates.
At each step, add the first item $v$ from $L$
into the independent set and delete from $L$ both $v$ and the nodes
which are adjacent to $v$.
The algorithm terminates when $L$ becomes empty.
Clearly, the running time of this implementation is $O(n \log{n} + m)$
where $n$ and $m$ are respectively the number of vertices and the
number of edges in the graph.

\subsection{Domination Problems}

In this section we give efficient heuristics for various domination
problems for unit disk graphs. 
The problems considered in this
section are minimum dominating set,
minimum independent dominating set,
minimum connected dominating set, and minimum total dominating set.

\subsubsection{Dominating Set and Independent Dominating Set}

An independent set $S$ is {\bf maximal} if no proper superset of $S$
is also an independent set.
It is well known and easy to see that for any graph, 
any maximal independent set is
also a dominating set.
Thus, one method of finding a dominating set (which is also an
independent dominating set) in a graph is to find a maximal
independent set.
A straightforward method of finding a maximal independent set
is to repeat the following steps until the graph becomes empty:
Select an arbitrary node $v$, add $v$ to the independent set
and delete $v$ and $N(v)$ from the graph.
Obviously, this method can be implemented to run in linear time.

We now observe that for unit disk graphs, the size of any
maximal independent set is within a factor of 5 of the size
of a minimum dominating set. 

\begin{theorem}\label{th:perf_adom}
Let $G$ be a unit disk graph.
Let $D^*$ be a minimum dominating set for $G$ and let 
$D$ be any maximal independent set for $G$.
Then $|D| \leq 5\,|D^*|$.
\end{theorem}

\noindent
{\bf Proof:} 
Since $D$ is an independent set, by Lemma~\ref{lem:k16},
no vertex in $D^*$ can dominate more than 5 vertices in $D$. 
Hence, $|D^*| \geq |D|/5$ and the theorem follows. \QED

Clearly, the same performance guarantee holds 
for independent domination because $D$ is also an independent set. 

In general, given a graph which is $K_{1,p}$ free, the cardinality
of any maximal independent set is at most 
$p-1$ times that of a minimum (independent) dominating set. 

\subsubsection{Connected Domination and Total Domination}

We assume that $G$ has no isolated nodes;
otherwise, $G$ has neither a connected dominating set nor a total 
dominating set.
If $G$ consists of two or more connected components each containing two or more
nodes, then $G$ does not have a connected dominating set.
However, in this case, we can find a near-optimal total dominating set
as follows.
We first find a maximal independent set $X$ for $G$.
Then, for each node $x \in X$, we choose a node $y \in V-X$ such that
$x$ is adjacent to $y$.
Clearly, this procedure leads to a total dominating set of size at most $2|X|$.
Hence from Theorem~\ref{th:perf_adom}, the size of the 
resulting total dominating set
is at most 10 times that of a minimum dominating set.

In view of the above, we assume for the remainder of this section
that $G$ is connected.
For any connected graph, a connected dominating set is also a total
dominating set.
Therefore, we focus on finding a connected dominating set.

\newcommand{\Begin}{\mbox{{\bf begin} }}
\newcommand{\Do}{\mbox{ {\bf do} }}
\newcommand{\End}{\mbox{{\bf end} }}
\newcommand{\For}{\mbox{{\bf for} }}
\newcommand{\Gets}{\mbox{ $\leftarrow$ }}
\newcommand{\To}{\mbox{ {\bf to} }}
\newcommand{\Output}{\mbox{{\bf output} }}

\newcommand{\DS}{\mbox{\em DS}}
\newcommand{\NS}{\mbox{\em NS}}

\begin{figure}
\oldspacing
\rule{6.5in}{0.01in}
\begin{tabbing}
123 \= 123 \= 123 \= 123 \= \kill
1. \> Arbitrarily pick a vertex $v \in V$. \\

2. \> Construct the breadth-first spanning tree $T$ of $G$ rooted at $v$. \\

3. \> Let $k$ be the depth of $T$. \\

4. \> Let $S_i$ denote the nodes at level $i$ in $T$, for all $0 \leq i \leq k$.
 \\

5. \> Set $\IS_0 = \{v\}$; $\NS_0 = \emptyset$. \\

6. \> \For $i=1$ \To $k$ \Do \Begin \\

7. \> \> $\DS_i = \{ v | v \in S_i$ and $v$ is dominated by some
                         vertex in $\IS_{i-1} \}$. \\

8. \> \> Pick a maximal independent set $\IS_i$ in $G(S_i - \DS_i)$. \\

9. \> \> $\NS_i = \{ u\, |\, u$ is the parent (in $T$) of some $v \in \IS_i \}$.
                   (Note that $\NS_i \subseteq S_{i-1}$.) \\
10.\> \End \\

11.\> \Output\ $(\cup_{i=0}^k \IS_i) \cup (\cup_{i=0}^k \NS_i)$
      as the connected dominating set.
\end{tabbing}
\caption{Details of Heuristic CDOM}
\label{fig:heur_cdom}
\rule{6.5in}{0.01in}
\newspacing
\end{figure}

We begin with an overview of our heuristic
for connected domination.
The heuristic first selects
an arbitrary vertex $v$ from the graph and  
constructs the breadth-first spanning tree $T$ of $G$ rooted at $v$. 
For any node $v_j$, let $l(v_j)$ denote the 
number of edges in the path from $v$ to $v_j$ in $T$.  
We will refer to $l(v_j)$ as the {\bf level} of $v_j$.
Let $k$ denote the maximum level in $G$.
Thus, the tree $T$ partitions the vertices of $G$
into sets $S_i$, $0 \leq i \leq k$, where  
$S_i$ is the set of nodes at level $i$.
The connected dominating set produced by the heuristic is the union
of two sets of nodes.
The first set is a maximal independent set for $G$ obtained
by appropriately selecting an 
independent set $\IS_i$ from each graph
$G(S_i)$, $1 \leq i \leq k$.
The second set of nodes is used to ensure connectivity.
The details of this heuristic, which we call CDOM,
appear in Figure~\ref{fig:heur_cdom}.
We are now ready to prove the correctness of the heuristic and
its performance guarantee.

\begin{lemma}
  The output of Heuristic CDOM is a connected dominating set.
\end{lemma}

\noindent
{\bf Proof:}
  We prove by induction on $k$ that $\cup_{i=0}^k \IS_i$ dominates
  the graph $G(\cup_{i=0}^k S_i)$, and that
  $G((\cup_{i=0}^k \IS_i) \cup (\cup_{i=0}^k \NS_i))$ is connected.

  {\bf Basis:} For $k=0$ we have $G(\cup_{i=0}^k S_i) = (\{v\}, \emptyset)$.
  Furthermore, $(\cup_{i=0}^k \IS_i) \cup (\cup_{i=0}^k \NS_i) 
  = \IS_0 \cup \NS_0 = \{v\}$,
  so the lemma holds.
 
  {\bf Inductive Step:} Assume the statment is true for some 
  $k \geq 0$ and consider $k+1$.
  By the inductive hypothesis, $\cup_{i=0}^k \IS_i$, 
  and therefore $\cup_{i=0}^{k+1} \IS_i$,
  dominates $G(\cup_{i=0}^k S_i)$.
  We therefore only need to show that $S_{i+1}$ is also dominated.
  By Step~7, all vertices in $\DS_{k+1}$ are dominated by $\IS_k$, and
  by Step~8, all vertices in $S_{k+1} - \DS_{k+1}$ 
  are dominated by $\IS_{k+1}$
  (since $IS_{k+1}$ is a maximal independent set in 
  $G(S_{k+1} - \DS_{k+1}))$.
  That is, $S_{k+1}$ is dominated by $\IS_k \cup \IS_{k+1}$, and therefore
  by $\cup_{i=0}^{k+1} \IS_i$.

  Now, by the inductive hypothesis,
  $G((\cup_{i=0}^k \IS_i) \cup (\cup_{i=0}^k \NS_i))$ is connected.
  By Step~9, $\NS_{k+1} \subseteq S_k$ since $\IS_{k+1} \subseteq S_{k+1}$.
  By the inductive hypothesis, $S_k$, and therefore $\NS_{k+1}$, is dominated
  by $(\cup_{i=0}^k \IS_i)$.
  Therefore $G((\cup_{i=0}^k \IS_i) \cup (\cup_{i=0}^{k+1} \NS_i))$ 
  is connected.
  Finally, by Step~9, every vertex in $\IS_{k+1}$ is adjacent to $\NS_{k+1}$.
  Therefore $G((\cup_{i=0}^{k+1} \IS_i) \cup (\cup_{i=0}^{k+1} \NS_i))$ 
  is connected. \QED

\begin{theorem}\label{th:perf_cdom}
Given a unit disk graph $G$, Heuristic CDOM computes
a connected dominating set whose size is no
more than 10 times that of an optimal connected
dominating set for $G$.
The same heuristic also yields a total dominating set of size no more than
10 times that of a minimum total dominating set.
\end{theorem}

\noindent
{\bf Proof:}
 Let $\IS = \cup_{i=0}^k \IS_i$ and $\NS = \cup_{i=0}^k \NS_i$.
  It is easy to verify that $\IS$ is a maximal independent set.
  Therefore, by Theorem~4.8, $|\IS|$ is no more than 5 times the size
  of a minimum dominating set (\OPT), and hence also within 5 times
  the size of minimum connected dominating set.
  By Steps~5 and 9, $|\NS_i| \leq |\IS_i|$ for all $0 \leq i \leq k$,
  so $|\NS \cup \IS| \leq 2|\IS| \leq 10|\OPT|$, and the bound of
  10 follows.
  Since every connected dominating set is also a total dominating set,
  the bound with respect to a minimum total dominating set is also 10.
\QED

In general, given a $K_{1,p}$ free graph
(i.e., a graph which has no subgraph isomorphic to $K_{1,p}$), the above
heuristic obtains a solution which is no more than $2(p-1)$ times any optimal
connected/total dominating set.  

\section{Extensions to other intersection graphs}

Thus far, we have considered intersection graphs of circles with the 
the same radius.
In this section, we show that the ideas presented in the previous
sections can be extended to obtain constant performance guarantees
for several optimization problems for intersection graphs of
circles of arbitrary radii in the plane and intersection graphs of
regular polygons. 

\subsection{Circle intersection graphs}

We use the term {\bf circle intersection graphs}
to refer to intersection graphs of circles of arbitrary radii
in the plane.
A graph $G$ is a circle intersection graph if the vertices
of $G$ can be placed in one-to-one correspondence with a set
of circles in the plane such that two vertices of $G$ are
adjacent if and only if their corresponding circles intersect.
As before, tangential circles are assumed to intersect.
We note that the class of circle intersection graphs is
different from the class of {\bf coin graphs} considered in
\cite{Th91}.
This is because in the case of coin graphs, two vertices
are adjacent if and only if their corresponding circles {\em touch}.
It is known that the class of coin graphs coincides with the
class of planar graphs \cite{Th91}.

We now show how to extend the heuristics for unit disk graphs to
circle intersection graphs.
We begin with a lemma which is similar to Lemma~\ref{lem:k16}.

\begin{lemma}\label{lem:circ_k16}
Let $G$ be a circle intersection graph.
Then, $G$ has a node $v$ such that the graph $G(\{v\} \cup N(v))$
does not contain an induced subgraph isomorphic to $K_{1,6}$.
\end{lemma}

\noindent
{\bf Proof:} Consider a circle intersection model for $G$.
Let $v$ be a node whose corresponding circle $C$ has the smallest
radius.
We can now use the argument given in the proof of Lemma~\ref{lem:k16}
to obtain the result. \QED

The above lemma indicates that every circle intersection graph $G$
contains a node $v$ such that the size of any independent set in
$G(N(v))$ is at most 5.
Using this observation and a proof similar to that of 
Lemma~\ref{lem:triangle_free}, it is easy to prove the following
result.

\begin{lemma}\label{lem:circ_triangle_free}
Any triangle-free circle intersection graph 
can be colored using 6 colors. \QED
\end{lemma}

Lemma~\ref{lem:circ_triangle_free} enables us to obtain a heuristic
for vertex cover with a performance guarantee of 5/3 for circle intersection
graphs.
The heuristic is the same as Heuristic VCover in Figure~\ref{fig:heur_vc}
except that in Step 6, the graph $G(Q)$ is colored 
using 6 colors (instead of 4).
The performance guarantee of 5/3 follows from Theorem~\ref{th:hoch} and
Lemma~\ref{lem:local}.

The heuristic given for off-line coloring of unit disk graphs (Section 4.2.1)
can be easily 
shown to provide a performance guarantee of 6 for circle
intersection graphs.
The proof is very similar to that of Theorem~\ref{th:perf_off_color}.
The only difference is the following.
We prove that $OPT(G) \geq \delta(G)/6 + 1$ by considering
the vertex $v^*$ to be a vertex whose corresponding circle has the smallest
radius and using Lemma~\ref{lem:circ_k16}.

It is also possible to obtain a heuristic with a performance
guarantee of 5 for the maximum independent set problem for 
circle intersection graphs.
This heuristic is identical to that for unit disk graphs (Section 4.3)
except that at each step, we select a node $v$ such that no independent set 
in the subgraph $G(N(v))$ is of size 6 or more.
Lemma~\ref{lem:circ_k16} enables us to conclude that the heuristic
provides a performance guarantee of 5.

The following theorem summarizes the above discussion.

\begin{theorem}\label{perf_circ_all}
There are polynomial time heuristics
that provide performance guarantees of 5/3, 6, and 5 respectively 
for minimum vertex cover,
minimum off-line vertex coloring, 
and maximum independent set problems for
circle intersection graphs. \QED
\end{theorem}

\subsection{Intersection graphs of regular polygons}

The ideas in the previous sections can also be extended to intersection
graphs of regular polygons. 
Let us call a $p$ sided polygon 
{\bf a unit regular polygon} if the 
polygon is inscribed in a circle of radius 1 
and the sides of the polygon are all equal.
Each such polygon can  be uniquely
specified up to rotation by the number of sides and the
coordinates of the center of the polygon. 
Now, for any graph which is an intersection graph of $p$-sided unit 
regular polygons, 
we can derive a (crude) bound on the size of a maximum independent
set in the subgraph induced on the neighborhood of any node. 
This bound
depends only on the number of sides $p$ of the polygon.

\begin{lemma}\label{lem:polygons}
Let $P$ be a $p$-sided unit regular polygon,
and let $S$ be a set of $p$-sided unit regular polygons such that 
every polygon in $S$ intersects $P$ and no two polygons in $S$
intersect each other.
There is an integer $t$, depending only on $p$, such that
$|S| \leq t$.
\end{lemma}

\noindent
{\bf Proof:} 
It is easy to verify that the area $\alpha_p$ of a $p$-sided unit regular
polygon is given by $\alpha_p = p \sin{(2\pi/p)}/2$.
It is also easy to see that any $p$-sided unit regular polygon which
intersects $P$ must lie completely within the circle $C$ of radius
3 with center at the center of $P$.
The area of this circle is $9\pi$.
Thus, the maximum number of $p$-sided unit regular polygons which
lie completely within $C$ and which do not intersect each other is
$\lceil\,{\rm Area~of~} C / \alpha_p\,\rceil$ which is
equal to $\lceil\,18\pi / (p \sin{(2\pi/p)})\,\rceil$.
Obviously, this bound depends only on $p$.
\QED

Thus, for any {\em fixed} $p$, there is a fixed integer $t$ such that
the intersection graphs of $p$-sided unit regular polygons 
do not contain $K_{1,t+1}$ as an induced subgraph.
Therefore, we can obtain polynomial time heuristics
with constant performance guarantees
for most of the problems considered 
in Section 4 for intersection graphs of
$p$-sided unit regular polygons. 
The following theorem provides a formal statement of our
results for intersection graphs of $p$-sided unit regular polygons. 

\begin{theorem}
For intersection graphs of $p$-sided unit regular polygons, there
are heuristics with constant performance guarantees
for the following optimization problems:
minimum vertex cover,  maximum independent set, 
minimum off-line coloring,
minimum dominating set, minimum independent dominating set,
minimum connected dominating set and minimum total dominating set. \QED
\end{theorem}

In a similar fashion, we can also extend the results in 
the previous sections to obtain 
heuristics with constant performance guarantees 
for intersection graphs of unit balls in higher dimensions. 

\section{Concluding remarks}

We have given efficient approximations for a variety of
standard problems on unit disk graphs. 
These heuristics have been designed
using the underlying geometric structure of unit disk graphs.
We observed that our heuristics can be extended 
to intersection graphs of circles in the plane and to 
intersection graphs of regular
polygons.

The heuristics presented in this paper do not require 
a geometric representation of the input graph.
In \cite{HMR+94}, we have shown that when a geometric representation
of a unit disk graph is available as part of the input, it is possible
to obtain polynomial time approximation schemes for maximum independent
set, minimum vertex cover and minimum dominating set problems.
These results are obtained using ideas from \cite{Ba83} and \cite{HM85}.

Our heuristics for the maximum independent set problem for unit disk graphs
and circle intersection graphs also provide a slight generalization of
a result in \cite{Ho83}.
In that paper, Hochbaum observed that for any graph which does not
contain $K_{1,r}$ as an induced subgraph, there is a heuristic with
a performance guarantee of $r-1$ for the maximum independent set 
problem.
Our results point out that if a graph $G$ has a node $v$ such that
$G(\{v\} \cup N(v))$ does not contain $K_{1,r}$ as an induced subgraph
and this property holds for all induced subgraphs of $G$,
then there is a heuristic with a performance guarantee of $r-1$
for the maximum independent set problem for $G$.

\vspace*{0.5in}

\noindent
{\bf Acknowledgements:}
We thank R. Ravi for pointing to us the work of Clark, Colbourn and Johnson 
\cite{CCJ90} and also for numerous suggestions during the course of
writing this paper. 
We would also like to thank
Professor Richard Stearns and 
Venkatesh Radhakrishnan for constructive suggestions.

\oldspacing

\baselineskip = \normalbaselineskip




\end{document}